\documentclass[11pt]{amsart} 

\usepackage[utf8]{inputenc} 

\usepackage[margin=1in]{geometry} 

\usepackage{graphicx,color} 

\allowdisplaybreaks[1]

\usepackage{pxfonts}
\usepackage{hyperref}

\usepackage{array} 
\usepackage{paralist} 
\usepackage{verbatim} 
\usepackage{subfig} 
\usepackage{amsfonts, amssymb, amsmath} 
\usepackage{amsthm}
\usepackage{hyperref} 
\usepackage{xcolor} 

\numberwithin{equation}{section} 

\newtheorem{theorem}{Theorem}[section]
\newtheorem{corollary}{Corollary}[theorem]
\newtheorem{lemma}[theorem]{Lemma}
\newtheorem{proposition}[theorem]{Proposition}

\title{Hölder regularity of solutions of degenerate parabolic equations of general dimension}

\author{Hyo Seok Jang}
\address{
The Research Institute of Basic Sciences, Seoul National University, Gwanak-ro 1, Gwanak-Gu, Seoul, 08826, South Korea 
}
\email{ hyoseok.jang@snu.ac.kr }

\author{Ki-Ahm Lee }
\address{
Department of Mathematical Sciences, Seoul National University, Gwanak-ro 1, Gwanak-Gu, Seoul, 08826, South Korea \&
Research Institute of Mathematics, Seoul National University, Gwanak-ro 1, Gwanak-Gu, Seoul, 08826, South Korea
}
\email{ kiahm@snu.ac.kr }

\begin{document}
\maketitle
\begin{abstract}
We establish the Alexandroff-Bakelman-Pucci estimate, the Harnack inequality, the Hölder regularity and the Schauder estimates to a class of degenerate parabolic equations of non-divergence form in all dimensions
\begin{equation} 
\mathcal{L}u:= u_t -Lu= u_t -(x a_{11} u_{xx} +2\sqrt{x} \sum_{j=2}^n a_{1j} u_{x y_j} + \sum_{i,j=2}^n a_{ij} u_{y_i y_j} + b_1 u_x +\sum_{j=2}^n b_j u_{y_j} ) =g\ 
\end{equation}
on \(x \geq 0, y=(y_2,\ldots, y_n) \in \mathbb{R}^{n-1}\), with bounded measurable coefficients.


\vspace{5mm} 
\end{abstract}

\section{Introduction}
We study the regularity of a family of degenerate parabolic equations of non-divergence form 

\begin{equation} \label{eq:degenerate:parabolic}
\mathcal{L}u:= u_t -Lu= u_t -(x a_{11} u_{xx} +2\sqrt{x} \sum_{j=2}^n a_{1j} u_{x y_j} + \sum_{i,j=2}^n a_{ij} u_{y_i y_j} + b_1 u_x +\sum_{j=2}^n b_j u_{y_j} ) =g
\end{equation}
on \(x\geq 0, y=(y_2,\ldots, y_n) \), with bounded and measurable coefficients with the conditions
\begin{equation} \label{assumption:coefficients}
a_{ij} \xi_i \xi_j \geq \lambda |\xi|^2 , \ |a_{ij}|, |b_i| \leq \lambda^{-1}, \ \text{ and }  \frac{2b_1}{a_{11}}\geq \nu>0
\end{equation}
for some constants \(0<\lambda<1, 0<\nu<1\), where the dimension \(n \geq 2\).

In this paper, we establish the Alexandrov-Bakelman-Pucci estimate, the Harnack inequality and the Hölder continuity of the equation \eqref{eq:degenerate:parabolic}, extending the precedant work in dimension two of Daskalopoulos and Lee \cite{daskalopoulos-lee03} to general dimensions, which is generalizing the result of Krylov and Safonov \cite{krylov-safonov80} and Tso \cite{tso85}. Because of the degeneracy of the equation \eqref{eq:degenerate:parabolic}, we use a singular metric \(s\) to scale all the estimates, as done in \cite{daskalopoulos-lee03}.

Also, we obtain Schauder estimates of the equation \eqref{eq:degenerate:parabolic} in general dimensions no less than two, which Daskalopoulos and Hamilton \cite{daskalopoulos-hamilton99} acquired in the case of dimension two.

Daskalopoulos and Lee \cite{daskalopoulos-lee03} showed the Hölder continuity of the solutions to the equation \eqref{eq:degenerate:parabolic} in dimension two. Using these results of \cite{daskalopoulos-lee03} for the model equation under certain coordinates, they in \cite{daskalopoulos-lee04} showed that the solution of the two-dimensional Gauss curvature flow of Firey \cite{firey74} with a flat side, which was introduced by Hamilton \cite{hamilton94}, exists smoothly and the interface is smooth for all time until the flat side vanishes. We can apply the results in this paper to obtain similar results for the Gauss curvature flows in higher dimesions as well.

If a flow by Gauss curvature with a flat side is represented by a graph, then its equation can be written by the fully non-linear equation
\begin{equation}
f_t =\frac{\det(D^2 f)}{(1+|Df|^2)^{(n+1)/2}}.
\end{equation}

As in the case of dimension two \cite{daskalopoulos-lee03} \cite{daskalopoulos-lee04}, if we assume that \(g=\sqrt{2f}\) vanishes linearly at the interface initially, then a solution which is either smooth or \(C_s^{2,\beta}\) for some \(0<\beta<1\) up to the interface, depending on the dimension, exists for all time until the flat side vanishes. We will introduce this result soon. We can show this by transforming the equation at a point on the interface between the flat side and the graph by 
\( \big(x_1, x_2 \ldots, x_n, g(x_1,x_2, \ldots, x_n, t)\big) \) to \( \big(h(x_{n+1}, x_2, \ldots, x_n, t), x_2, \ldots, x_n, x_{n+1} \big)\). With the notation \( \mathcal{I}=h_{n+1}^2 +x_{n+1}^{2\alpha } +x_{n+1}^{2\alpha } \sum_{i=2}^n h_i^2 =h_{n+1}^2  I =h_{n+1}^2  (1+|Df|^2)\), the evolution of \(h\) is
\begin{equation} \label{eq:h_t}
\begin{split}
h_t =& -h_{n+1}^{n+2}  \frac{x_{n+1}^{-(n-1)} }{ \mathcal{I}^{(n+1)/2} } \det \bigg( -x_{n+1}\frac{1}{h_{n+1} }\Big( \frac{h_i h_j}{h_{n+1}^2}h_{n+1, n+1}-\frac{h_i}{h_{n+1} } h_{n+1, j}-\frac{h_j}{h_{n+1} } h_{n+1, i}+h_{ij} \Big)+\alpha \frac{h_i h_j }{h_{n+1}^2 }\bigg)
\end{split}
\end{equation}
and the linearized equation of \eqref{eq:h_t} satisfies a degenerate equation of type \eqref{eq:degenerate:parabolic} under suitable geometric conditions.

There have been a number of notable results for similar types of degenerate equations. Local a priori \(C_s^{2,\alpha}\)-estimates for degenerate equations of the form
\begin{equation} \label{eq:deg:hom}
u_t-\big( x( a_{11} u_{xx} +2a_{12} u_{xy} +a_{22} u_{yy} )+b_1 u_x+b_2 u_y \big) =g
\end{equation}
with \(C_s^{\alpha}\) coefficients with the elipticity condition \(a_{ij} \xi_i \xi_j \geq \lambda |\xi|^2\) and a lower bound constant \(\nu>0\) for which \(b_1 \geq \nu\) have been shown by Daskalopoulos and Hamilton \cite{daskalopoulos-hamilton88}, as a main step on establishing the short-time existence of a smooth-up-to-the-interface solution \cite{daskalopoulos-hamilton88} of the porous medium equation \eqref{eq:pme} 
\begin{equation} \label{eq:pme}
f_t=f\Delta f +\nu|Df|^2, \nu>0
\end{equation}
with suitable \(C_s^{2,\alpha}\) data,  which describes the pressure \(f\) of a gas through a porous medium. Daskalopoulos, Hamilton and K.-A. Lee \cite{daskalopoulos-hamilton-lee01} also showed the all-time \(C^{\infty}\) regularity of solutions to \eqref{eq:pme} using the Hölder a priori estimates of Koch \cite{koch99} of solutions to degenerate equations of the divergence form
\begin{equation}
x_n \Delta_{\mathbb{R}^{n-1}}  u-x_n^{-\sigma} \partial_{x_n} (x_n^{1+\sigma} \sum_{j=1}^n a^j\partial_j u) -u_t =g
\end{equation}
where Koch showed such estimates via scaling a Moser's iterattion argument according to a singular metric \(s\), chosen suitably for the problem. It is also worth mentioning that Lin and Wang \cite{lin-wang98} have established \(C^{\alpha}, C^{1,\alpha} \) and \(C^{2,\alpha}\) regularity of solutions to degenerate elliptic equations of the type \eqref{eq:deg:hom} with \(b_1\leq 0\).

In higher dimensions, we have \(C_s^{2,\alpha}\) regularity of the degenerate equations
\begin{equation}
u_t-x_n^\gamma \Delta u=g
\end{equation}
with a constant \(0<\gamma<1\) from Kim and Lee \cite{kim-lee09}. 

Recently, T. Kim, K.-A. Lee and H. Yun \cite{kim-lee-yun23} obtained the boundary \(C^{2,\alpha}_s\) regularity and higher regularity for degenerate equations of the form
\begin{equation}
u_t-(x^{1+\delta} a_{11} u_{xx} +2x^{\frac{1+\delta}{2}} \sum_{j=2}^n a_{1j} u_{x y_j} + \sum_{i,j=2}^n a_{ij} u_{y_i y_j} + x_1^{\frac{1+\delta}{2}} b_1 u_x +\sum_{j=2}^n b_j u_{y_j} +cu) =g
\end{equation}
with a constant \(0<\delta<1\), which is derived from the Constant Elasticity of Variance model in mathematical finance. They proved the smoothness of solution under certain smoothness conditions on coefficients and forcing terms as well.

We can ask later the regularity of various types of degenerate equations of the general form
\begin{equation}
u_t-\big(\sum_{i,j=1}^n x^{\alpha_i} x^{\alpha_j} a_{ij} u_{ij} +\sum_{i=1}^n b_i u_i +cu \big)=g
\end{equation}
with weaker conditions on coefficients, which arise in diverse fields of science and engineering.

\subsection{Main results} 

We introduce a priori estimates, which are the main results of this article. Let us denote by \(d\mu \) the measure \(d\mu =x^{\frac{\nu}{2}-1} dx dy_2 \ldots dy_n \) in the following statements.
\begin{theorem}[Hölder regularity] \label{theorem:hoelder:lin} 
Assume that the coefficients of the operator \( L\) are smooth on \(\mathcal{C}_\rho, \ \rho>0 \), and satisfy the bound conditions. Then there exists a constant \( 0 <\alpha<1
\), so that for any \(r<\rho\)
\begin{equation*}
\Vert u \Vert_{C_\rho^\alpha (\mathcal{C}_r) } \leq C(n,r,\rho) \left( \Vert u \Vert_{C^0 (\mathcal{C}_1) } + \left( \int_{\mathcal{C}_\rho} g^{n+1}(x,t) d\mu dt\right)^{1/(n+1)}\right)
\end{equation*}
for all smooth functions \(u\) on \(\mathcal{C}_\rho \) for which \(Lu-u_t=g\).
\end{theorem}

\begin{theorem} [Schauder Estimate] \label{theorem:schauder}
Let \( k\) be a nonnegative integer and let \(0<\alpha<1 \) and \(v>0\) be given. Then, for any \( r<1\) there exists a constant \(C >0 \) depending on \( k , \alpha, v\) and \(r \) such that 
\begin{align*}
\Vert u\Vert_{\mathcal{C}_s^{k,2+\alpha} (\mathcal{B}_r)} \leq C( \ \Vert u\Vert_{\mathcal{C}_s^{0} (\mathcal{B}_1)} +\Vert L_0 u\Vert_{\mathcal{C}_s^{k,\alpha} (\mathcal{B}_1)} \ )
\end{align*}
for any smooth function \(u\) on \( \mathcal{B}_1 \).
\end{theorem}
  
\subsection{Summary}

The outline of this paper is as follows. We show the Hölder a priori regularity Theorem \ref{theorem:hoelder:lin} of the equation \eqref{eq:degenerate:parabolic} by showing the Alexandrov-Bakelman-Pucci estimate Lemma \ref{lem:ABP:lin} and the Harnack inequality Lemma \ref{lemma:harnack}. We also derive Schauder estimates Theorem \ref{theorem:schauder} and the short-time existence Theorem \ref{theorem:short:variable} of the model equation \eqref{eq:degenerate:parabolic}.

\subsection{Notations}
Here are some notations which we will use throughout the paper:
\begin{itemize}
\item \( \mathcal{B}_\eta=\mathcal{B}_\eta(P) :=\{ (x, (y_2, \ldots, y_n ),t) =(x,y,t)| x \geq 0, | x -x_0 | \leq \eta^2, |y-y_0|\leq \eta, t_0-\eta^2 \leq t \leq t_0 \} \) near \( P=(x_0,y_0,t_0)\) for \( \eta>0\). \( \mathcal{B}_r^\gamma =\mathcal{B}_r^\gamma (Q):=\{ x \geq 0, | \sqrt{x} -\sqrt{x_0}| \leq r, \ \gamma |y-y_0^\prime|\leq r, t_0^\prime -r^2 \leq t \leq t_0^\prime \} \) near \( Q=(z_0,y_0^\prime,t_0^\prime)\) for \( \gamma>0, r>0\). \( \mathcal{C}_\rho:=\{ x \geq 0, |x-x_0|^2 \leq \rho^2, |y-y_0|^2 \leq \rho^2, t_0-\rho^2 \leq t  \leq t_0 \} \), \( \mathcal{Q}_\rho:=\{ x\geq 0, |\sqrt{x}-\sqrt{x_0}| \leq \rho, |y-y_0|^2 \leq \rho^2, t_0-\rho^2 \leq t  \leq t_0 \} \) for \( \rho>0\). 
\end{itemize}

\section{Hölder regularity} \label{sec:holder}
Throughout the section, we assume that \( u\) is a solution of the equation \eqref{eq:degenerate:parabolic} and it is smooth in time \( 0  \leq t \leq T < T_c \).

\subsection{Alexandrov-Bakelman-Pucci estimate}

The proof of Theorem \ref{theorem:hoelder:lin} needs the following \\ Alexandrov-Bakelman-Pucci (A-B-P) estimate, which is a generalization of Theorem 3.2 in \cite{daskalopoulos-lee03} in dimension \(n\).

\begin{lemma}[Alexandrov-Bakelman-Pucci estimate] \label{lem:ABP:lin}
Let \(u\) be a smooth function on \(\mathcal{C}_\rho \) of the equation \(Lu-u_t=g\) with coefficients of the operator \( L\) satisfying the bound conditions. Additianlly, assume that \(u\leq 0\) on \( \{ |s-s_0|=\rho, |y-y_0|=\rho, t-t_0=\rho^2, \text{ where } s,t \in \mathbb{R}, y\in \mathbb{R}^{n-1} \} \cap \mathcal{C}_\rho\). Then there exists a constant \( 0 <\alpha<1
\), so that for any \(r<\rho\)
\begin{equation} \label{eq:ABP:lin}
\sup_{\mathcal{C}_\rho} u^+  \leq C(n,\lambda,\nu) \rho^{n/(n+1)} \rho_\nu^{1/(n +1)} \left( \int_{\Gamma^-} (g^-)^{n+1}(s,y,t) s^{\nu-1} ds dy dt \right)^{1/(n+1)}
\end{equation}
with
\begin{equation} \label{nota:rho-nu}
\rho_\nu(s_0)=(s_0+\rho)^{2-\nu}-s_0^{2-\nu}.
\end{equation}
\end{lemma}

\begin{proof}
The proof refers to that of Theorem 3.2 in \cite{daskalopoulos-lee03} in dimension \(n\). We consider the variable \( z =\frac{s^{2-\nu}}{2-\nu}\) and the gradient map
\begin{equation}
Z(z,y,t)=\big(u_z, D_y u, u -(zu_z +y \cdot D_y u) \big)
\end{equation}

with
\begin{align*}
\frac{\partial Z}{\partial(z,y,t)}
& =\begin{pmatrix}
& u_{zz} \ \ \ \ & D_y u_z \ \ \ \  & -( z u_{zz} +y \cdot D_{y} u_z) \\ 
& D_{y} u_z \ \ \ \ & D_y^2 u \ \ \ \  & -(z D_y u_{z} +y \cdot D_y (D_y u) ) \\
& u_{zt} \ \ \ \ & D_y u_t \ \ \ \ & u_t -(zu_{zt} +y \cdot D_y u_t)
\end{pmatrix}  
\end{align*}
where \( -(z D_y u_{z} +y \cdot D_y (D_y u) )\) is the \( (n-1 ) \times 1 \) matrix whose \(i\)-th entry is \( -(z u_{y_i z} +y \cdot D_y u_{y_i}  ) \).

So we have
\begin{align*}
\det \frac{\partial Z}{\partial(z,y,t)}
& = (-1)^{n-1} u_t \det \begin{pmatrix}
& u_{zz} \ \ \ \ & D_y u_z \ \ \ \  \\ 
& D_{y} u_z \ \ \ \ & D_y^2 u \ \ \ \ 
\end{pmatrix} 
 +(-1)^{n} (zu_{zt} +y \cdot D_y u_t)\det \begin{pmatrix}
& u_{zz} \ \ \ \ & D_y u_z \ \ \ \  \\ 
& D_{y} u_z \ \ \ \ & D_y^2 u \ \ \ \ 
\end{pmatrix}  \\
&  -\sum_{i=1}^{n-1} (-1)^{i} (z u_{y_i z} +y \cdot D_y u_{y_i} )
\det \begin{pmatrix}
& u_{zz} \ \ \ \ & D_y u_z \ \ \ \ \\ 
& [ D_{y} u_z \ \ \ \ & D_y^2 u ]_{-i} \ \ \ \ \\
& u_{zt} \ \ \ \ & D_y u_t \ \ \ \
\end{pmatrix}  \\
& -( z u_{zz} +y \cdot D_{y} u_z) 
\det \begin{pmatrix}
& D_{y} u_z \ \ \ \ & D_y^2 u \ \ \ \ \\
& u_{zt} \ \ \ \ & D_y u_t \ \ \ \
\end{pmatrix}  \\
& = (-1)^{n-1} u_t \det \begin{pmatrix}
& u_{zz} \ \ \ \ & D_y u_z \ \ \ \  \\ 
& D_{y} u_z \ \ \ \ & D_y^2 u \ \ \ \ 
\end{pmatrix}  
 +\det \begin{pmatrix}
& u_{zz} \ \ \ \ & D_y u_z \ \ \ \  & -( z u_{zz} +y \cdot D_{y} u_z) \\ 
& D_{y} u_z \ \ \ \ & D_y^2 u \ \ \ \  & -(z D_y u_{z} +y \cdot D_y (D_y u) ) \\
& u_{zt} \ \ \ \ & D_y u_t \ \ \ \ & -(zu_{zt} +y \cdot D_y u_t)
\end{pmatrix} 
\end{align*}
where \( [ D_{y} u_z \ \ \ \ D_y^2 u ]_{-i} \) is the \( (n-2 ) \times n \) matrix we get by removing the \(i\)-th row \( ( u_{y_i z} \ \ \ \ D_y u_{y_i} ) \) from the \( n \times n \) matrix \(  ( D_{y} u_z \ \ \ \ D_y^2 u  ) \). Because the \(i\)-th entry \( -y \cdot D_y u_{y_i}  = -\sum_{j=1}^{n-1} y_j u_{y_i y_j}  \) of \(y \cdot D_y u_{y_i}\) is \(y_j \) times the \((i,j)\)-th entry of \(D_y^2 u\), we have 
\begin{align*}
\det \begin{pmatrix}
& u_{zz} \ \ \ \ & D_y u_z \ \ \ \  & -( z u_{zz} +y \cdot D_{y} u_z) \\ 
& D_{y} u_z \ \ \ \ & D_y^2 u \ \ \ \  & -(z D_y u_{z} +y \cdot D_y (D_y u) ) \\
& u_{zt} \ \ \ \ & D_y u_t \ \ \ \ & -(zu_{zt} +y \cdot D_y u_t)
\end{pmatrix} =0
\end{align*}
and hence
\begin{equation}
\det \left(\frac{\partial Z}{\partial(z,y,t)} \right)
=(-1)^{n-1} u_t \det \left(\frac{\partial (u_z, u_y)}{\partial(z,y)}\right).
\end{equation}

We define the sets
\begin{equation} \label{nota:gamma}
\begin{split}
\Gamma^+ = & \{ (s,y,t)\in \mathcal{C}_\rho: \ \frac{\partial (u_z, u_y)}{\partial(z,y)} \leq 0, \ u_z \leq 0, u_t \geq 0 \} \\
\Gamma^- = & \{ (s,y,t)\in \mathcal{C}_\rho: \ \frac{\partial (u_z, u_y)}{\partial(z,y)} \geq 0, \ u_z \geq 0, u_t \geq 0 \} 
\end{split}
\end{equation}
and we denote by \(\mathcal{C}_r(s_0, y_0, t_0) \) the cube
\begin{equation}
\mathcal{C}_r (s_0, y_0, t_0) =\{(s,y,y)\in \mathbb{R} \times  \mathbb{R}^{n-1} \times \mathbb{R} : \ s \geq 0, |s-s_0|\leq r, |y-y_0|\leq r, -r^2 \leq t \leq t_0\} 
\end{equation}
for any point \((s_0, y_0, t_0) \) with \(s_0\geq0 \) and any number \(r>0\).

\begin{equation}
E=\begin{pmatrix}
s^{2(1-\nu)}u_{zz} \ \ s^{1-\nu}u_{z y_2}  \ \ \cdots \ \ s^{1-\nu}u_{z y_n} \\
s^{1-\nu}u_{y_2 z} \ \ u_{y_2 y_2} \ \ \ \ \cdots \ \ u_{y_2 y_n}  \\
\ \ \ \vdots \ \ \ \ \ \ \ \ \ \ \ \ \vdots \ \ \ \ \ \ \ \ddots \ \ \ \ \ \ \vdots \\
s^{1-\nu}u_{y_n z} \ \ u_{y_n y_2} \ \ \ \ \cdots \ \ u_{y_n y_n} 
\end{pmatrix} 
\end{equation}

We have \( u_t \geq 0\) and \(D^2 u \leq 0\) on \(\Gamma^+\), so \(|u_t \det E|=u_t (-\det E)\). By the matrix formula (9.10) in Gilbarg-Trudinger
\begin{equation}
\det A \det B \leq \big( \frac{\text{trace} AB }{n}\big)^n, \ \ A, B\geq 0, \ \text{symmetric}, \ n \times n, 
\end{equation}
we have, for 
\begin{equation}
A=\begin{pmatrix}
& 1 \ \ \ \ \ \ \ \ & 0  & \cdots \ \ & 0 \\
& 0 \ \ \ \ \ \ \ & a_{11} \ \ \ \ & \cdots \ \ & a_{1 n}  \\
& \ \vdots \ \ \ \ \ \ \ \ & \vdots \ \ \ \ \ \ & \ddots \ \  \ \ & \vdots \\
& 0 \ \ \ \ \ \ \ & a_{n 1} \ \ \ \ & \cdots \ \ \ \ & a_{n n } 
\end{pmatrix} 
\end{equation} 
\begin{equation}
B =\begin{pmatrix}
& -u_t \ \ & 0 \ \ & 0  \ \ & \cdots \ \ & 0 \\
& 0 \ \  & s^{2(1-\nu)}u_{zz} \ \ & s^{1-\nu}u_{z y_2}  \ \ & \cdots \ \ & s^{1-\nu}u_{z y_n} \\
& 0 \ \ & s^{1-\nu}u_{y_2 z} \ \ & u_{y_2 y_2} \ \ \ \ & \cdots \ \ & u_{y_2 y_n}  \\
& 0 \ \ & \vdots \ \ \ \ \ \ \ \ \ \ \ \ & \vdots \ \ \ \ \ \ \ & \ddots \ \ \ \ & \vdots \\
& 0 \ \ & s^{1-\nu}u_{y_n z} \ \ & u_{y_n y_2} \ \ \ \ & \cdots \ \ & u_{y_n y_n} 
\end{pmatrix} 
\end{equation}

\begin{equation}
(n+1) \big(u_t\det(a_{ij})\cdot (-\det E ) \big)^{\frac{1}{n+1}} \leq \big( \text{trace} (a_{ij}) \cdot E  -u_t\big)
\end{equation}

Considering the relation \( z =\frac{s^{2-\nu}}{2-\nu}\), we know that \(u_s=s^{1-\nu} u_z\) and \( u_{ss}=(1-\nu)s^{-\nu} u_z+s^{2(1-\nu)} u_{zz} = \frac{(1-\nu)}{s} u_s+s^{2(1-\nu)} u_{zz}\) so that
\begin{equation}
E_{11}=s^{2(1-\nu)}u_{zz} =u_{ss}+\frac{(\nu-1)}{s} u_s, \  E_{1i}=E_{i1}=s^{1-\nu}u_{z y_i}=u_{s y_i}, \ E_{ij}=u_{y_i y_j} \ i,j =2, \ldots,n
\end{equation}
so
\begin{equation}
B =\begin{pmatrix}
& -u_t \ \ & 0 \ \ & 0  \ \ & \cdots \ \ & 0 \\
& 0 \ \  & u_{ss}+\frac{(\nu-1)}{s} u_s \ \ & u_{s y_2}  \ \ & \cdots \ \ & u_{s y_n} \\
& 0 \ \ & u_{y_2 s} \ \ & u_{y_2 y_2} \ \ \ \ & \cdots \ \ & u_{y_2 y_n}  \\
& 0 \ \ & \vdots \ \ \ \ \ \ \ \ \ \ \ \ & \vdots \ \ \ \ \ \ \ & \ddots \ \ \ \ & \vdots \\
& 0 \ \ & u_{y_n s} \ \ & u_{y_n y_2} \ \ \ \ & \cdots \ \ & u_{y_n y_n} 
\end{pmatrix} 
\end{equation}
and
\begin{equation}
\begin{split}
\text{trace} (a_{ij}) \cdot E  -u_t 
& =a_{11}\big( u_{ss}+\frac{(\nu-1)}{s} u_s\big) +\sum_{i=2}^n a_{1i} u_{y_i s } +\sum_{i=2}^n a_{i 1} u_{s y_i } +\sum_{j=2}^n\sum_{i=2}^n a_{ji} u_{y_i y_j } \\
& =a_{11}\big( u_{ss}+\frac{(\nu-1)}{s} u_s\big) +2\sum_{i=2}^n a_{1i} u_{y_i s } +\sum_{j=2}^n\sum_{i=2}^n a_{ji} u_{y_i y_j }
\end{split}
\end{equation}

Because \( \frac{b_1}{ 2a_{11} } \geq \nu >0 \) and \( u\) is a classical subsolution of equation
\begin{equation}
L_s u := a_{11} u_{ss} + 2\sum_{i=2}^n a_{1i} u_{y_i s } +\sum_{j=2}^n\sum_{i=2}^n a_{ji} u_{y_i y_j } +\frac{a_{11}}{s} \big( \frac{b_1}{ 2a_{11} } -1 \big) u_s +\sum_{i=2}^n b_i u_{y_i} \geq g
\end{equation}
we have
\begin{equation}
(n+1) \big(u_t\det(a_{ij})\cdot |\det E| \big)^{\frac{1}{n+1}} \leq g^{-1} +\sum_{i=2}^n |b_i| |u_{y_i} |
\end{equation}
and by Hölder's inequality
\begin{equation}
(n+1) \big(u_t\det(a_{ij})\cdot |\det E| \big)^{\frac{1}{n+1}} 
\leq ( k^{n+1} (g^{-1})^{n+1}  +\sum_{i=2}^n |b_i|^{n+1}  )^{\frac{1}{n+1} }
( k^{-\frac{n+1}{n} }  +\sum_{i=2}^n |u_{y_i} |^{\frac{n+1}{n} } )^{\frac{n}{n+1}  }
\end{equation}
for all numbers \( k>0\). Since \(  \det (a_{ij}) \geq \lambda^n \geq \lambda^{n+1}\), we see that
\begin{equation}
\big(u_t\cdot |\det E| \big)^{\frac{1}{n+1}} 
( k^{-\frac{n+1}{n} }  +\sum_{i=2}^n |u_{y_i} |^{\frac{n+1}{n} } )^{-\frac{n}{n+1}  }
\leq \frac{1}{n+1} \lambda^{-1} ( k^{n+1} (g^{-1})^{n+1}  +\sum_{i=2}^n |b_i|^{n+1}  )^{\frac{1}{n+1} }
\end{equation}

Let us define the function \(G \) on \(\mathbb{R}^{n+1} \) by
\begin{equation}
G(\xi,\zeta, \tau)=( k^{-\frac{n+1}{n} }  +|\xi |^{\frac{n+1}{n} } )^{-n  }
\end{equation}
for \( (\xi,\zeta, \tau) \in \mathbb{R}^{n-1} \times \mathbb{R} \times \mathbb{R} \). We assume that 
\begin{equation}
M=\max_{\mathcal{C}_\rho } u^+
\end{equation}
is finite.
Then 
\begin{equation}
D =\left[-\frac{cM}{\rho_\nu (s_0)},0 \right] \times \left[-\frac{cM}{\rho}, \frac{cM}{\rho} \right]^{n-1}  \times \left[-\frac{cM}{\rho}, \frac{cM}{\rho} \right] \subset Z(\Gamma^+)
\end{equation}
for some uniform constant \(c>0\).
We have
\begin{equation}
\int_D G 
\leq \int_{\Gamma^+ } G(Z) \left| \det \left( \frac{\partial Z}{\partial (s,y,t) } \right) \right| ds dy dt
= \int_{\Gamma^+ } ( k^{-\frac{n+1}{n} }  +|D_y u |^{\frac{n+1}{n} } )^{-n  }  \left| u_t \det E \right| d\mu dt
\end{equation}

Using the bound \( \sum_{i=2}^n |b_i|^{n+1} \leq \lambda^{-1} \), we see that
\begin{equation} 
\int_D G 
\leq \frac{1}{(n+1)^{n+1} \lambda^{n+1}}\int_{\Gamma^+ }  ( k^{n+1} (g^{-1})^{n+1}  +\lambda^{-(n+1) }  )  d\mu dt
\end{equation}
while
\begin{align*}
\int_D G = & \frac{c\rho}{\rho_\nu(s_0)} \int_{B_{\frac{cM}{\rho}} } ( k^{-\frac{n+1}{n} }  +|\xi |^{\frac{n+1}{n} } )^{-n  }d\zeta d\xi d\tau \\
\geq & \frac{c\rho}{\rho_\nu(s_0)} \int_{B_{\frac{cM}{\rho}} } ( k^{-(n+1) }  +|\xi |^{n+1 } + \zeta^{n+1 } +\tau^{n+1 })^{-1  }d\zeta d\xi d\tau 
=\frac{c\rho}{\rho_\nu(s_0)} \log 
\left( 1+\frac{c^{n+1} k^{n+1} M^{n+1} }{\rho^{n+1}} \right)
\end{align*}

Setting \( k\) by \( k^{-(n+1)} =\lambda^{n+1} \int_{\Gamma^+} (g^-)^{n+1} d\mu dt \), we can conclude that
\begin{equation}
\frac{c\rho}{\rho_\nu(s_0)} \log 
\left( 1+\frac{c^{n+1} k^{n+1} M^{n+1} }{\rho^{n+1}} \right) \leq C(n,\lambda, \nu).
\end{equation}

As \(\alpha=\frac{\rho}{\rho_\nu(s_0)} \geq 1\) for \( s_0<1, \rho<1\), it holds that
\begin{equation}
\log \left( 1+\frac{c^{n+1} k^{n+1} M^{n+1} }{\rho^{n}\rho_\nu(s_0)} \right) \leq C(n,\lambda, \nu),
\end{equation}
and consequently
\begin{equation}
M\leq C(n,\lambda, \nu) \rho^{\frac{n}{n+1}}\rho_\nu(s_0)^{\frac{1}{n+1}} \left( \int_{\Gamma^+} (g^-)^{n+1} d\mu dt \right)^{\frac{1}{n+1}}
\end{equation}
showing the desired result \eqref{eq:ABP:lin}.

\end{proof}

\subsection{Barrier construction}

We set the variable \(s=\sqrt{x}\) and define the operator \(\mathcal{L}_s u\) as
\begin{equation}
\mathcal{L}_s u =a_{11} u_{ss} +2 a_{1i} u_{s y_i} + a_{ij} u_{y_i y_j} +\frac{a_{11}}{s}\big[ \frac{b_1}{2a_{11}} u_x -1\big]+b_j u_{y_j} -u_t
\end{equation}
satisfying the conditions
\begin{equation}
|b_i|\leq \lambda^{-1} \ \text{ and } \ \frac{b_1}{2a_{11}}\geq \nu>0.
\end{equation}

We define a distance \( d_\gamma^2 ((x,y),(x_0,y_0)):=( \sqrt{x}-\sqrt{x_0})^2  +\gamma^2\sum_{i=2}^n (y_i-(y_0)_i)^2 \) between two points \((x,y),(x_0,y_0) \in \mathbb{R} \times \mathbb{R}^{n-1}\).
For any point \( (x_0, y_0) \in \mathbb{R} \times \mathbb{R}^{n-1} \) and any number \(0 <\rho \leq 1\), we denote \( K_{3\sqrt{2}\rho} =\mathcal{B}_{3\sqrt{2} \rho} \times (0,18\rho^2) \), \(\mathcal{Q}^1_{\frac{\rho}{2}} = \mathcal{Q}_{\frac{\rho}{2}}(x_0, y_0, \frac{\rho^2}{4} )\), \(\mathcal{Q}^2_{\frac{3\rho}{2}} = \mathcal{Q}_{\frac{3\rho}{2}}(x_0, y_0, \frac{10\rho^2}{4} )\).

We state the existence of a barrier function, the analogue of Lemma 3.4 in \cite{daskalopoulos-lee03} in dimension \(n\).
\begin{lemma} \label{lem:barrier:ineq}
For any point \( (x_0, y_0) \in \mathbb{R} \times \mathbb{R}^{n-1} \) with \( 0\leq x_0 \leq 1\) and any number \(0 <\rho \leq 1\), there exists a function \(\phi_\rho\) on  \(K_{3\sqrt{2}\rho} \), satisfying
\begin{equation} \label{ineq1:phi-rho}
\begin{cases}
\phi_\rho \geq 1 \text { in } \mathcal{Q}^2_{\frac{3\rho}{2}} \\
\phi_\rho \leq 0 \text { on } \partial_p K_{3\sqrt{2}\rho},
\end{cases}
\end{equation}

\begin{equation} \label{ineq2:phi-rho}
L\phi_\rho -\partial_t \phi_\rho \geq 0 \text { in } K_{3\sqrt{2}\rho} \setminus \mathcal{Q}^1_{\frac{\rho}{2}},
\end{equation}

and for a constant \( C(\lambda, \nu) >0 \)
\begin{equation}
\Vert \phi_\rho \Vert_{C^{1,1} (K_{3\sqrt{2}\rho})} \leq \frac{C(\lambda, \nu)}{\rho^2}.
\end{equation}

\end{lemma}

\begin{proof}

We show the case \(\rho=1\), and we will show the general case by a dialation at the end. Let us e define a distance \( \overline{d}_\gamma^2 ((x,y),(x_0,y_0)):=\frac{( x-x_0)^2 }{x+x_0} +\gamma^2\sum_{i=2}^n (y_i-(y_0)_i)^2 \) between two points \((x,y),(x_0,y_0) \in \mathbb{R} \times \mathbb{R}^{n-1}\). We fix the number \(\gamma>0 \) and simply denote \( \overline{d}_\gamma^2 ((x,y),(x_0,y_0))\) by \( \overline{d} ((x,y),(x_0,y_0)) \). We observe the two distance functions \(d_\gamma\) and \( \overline{d} \) are equivalent in that
\begin{equation} \label{ineq:dist:equiv}
d_\gamma \leq \overline{d} \leq \sqrt{2} d_\gamma 
\end{equation}

Given the base point \((x_0,y_0) \in \mathbb{R} \times \mathbb{R}^{n-1} \), let us consider the function \( \omega\) on points \((x,y, t) \in \mathbb{R} \times \mathbb{R}^{n-1} \times \mathbb{R}\)
\begin{equation}
\omega(x,y,t) =[18- \overline{d}^2((x,y),(x_0,y_0))] \Lambda(x,y,t)
\end{equation}
where
\begin{equation}
\Lambda(x,y,t) =\frac{1}{4\pi t} e^{-\frac{\overline{d}^2((x,y),(x_0,y_0))}{t}}.
\end{equation}

Let \(0<\tau_0<1, m>1\) and \(l>1\) be numbers to be determined later. Then we define a function
\begin{equation}
v(x,y,t)=e^{-mt} \omega^l(x,y,t+\tau_0)-M(\tau_0)
\end{equation}
where 
\begin{equation}
M(\tau_0)=\sup \{ \omega^l(x,y,\tau_0) : \overline{d}((x,y),(x_0,y_0))\geq \frac{1}{2}\}.
\end{equation}

Then we have \(v\leq 0\) on \(\partial_p K_{3\sqrt{2}} \setminus \mathcal{Q}^1_{\frac{1}{2}} \) by the equivalence relation between distances \eqref{ineq:dist:equiv}. Setting the number \(\tau_0>0\) sufficiently small, depending only on \(\gamma\), we can make \(v>0\) on \(\mathcal{Q}^2_{\frac{3}{2}}\).

From now on, we write \(\theta(x,y)=\overline{d}^2((x,y),(x_0,y_0))\) for the sake of simplicity of notation so that
\begin{equation}
\omega=(18-\theta)\Lambda, \ \text{ where } \Lambda =\frac{1}{4\pi t} e^{-\frac{\theta}{t}}.
\end{equation}

For the parabolic operator \(\mathcal{L}\) in the equation \eqref{eq:degenerate:parabolic}
\begin{equation}
\mathcal{L} v := v_t -Lv=v_t-(x a_{11} v_{xx} +2\sqrt{x} a_{1i} v_{x y_i} + a_{ij} v_{y_i y_j} + b_1 v_x +b_j v_{y_j} )
\end{equation}
we denote the coefficients \(\tilde{a}_{ij}\) by
\begin{equation}
\tilde{a}_{11}=x a_{11}, \  \text{ and } \tilde{a}_{1i}=\sqrt{x} a_{1i} , \ \tilde{a}_{ij}=a_{ij} \text{ for } i,j=2,\ldots,n,
\end{equation}
so that
\begin{equation}
\mathcal{L} v = v_t -Lv=v_t-(\sum_{i,j=1}^n \tilde{a}_{ij} v_{ij} + b_1 v_x +b_j v_{y_j} )
\end{equation}
where \(v_1=v_x\) and \(v_j=v_{y_j}\) for \(j =2, \ldots,n\).

Computing directly
\begin{equation}
\mathcal{L} v =\mathcal{L} (e^{-mt} \omega^l(x,y,t+\tau_0)-M(\tau_0) ),
\end{equation}
we see that
\begin{equation} \label{eq:omega}
\begin{split}
\mathcal{L} v = & \partial_t (e^{-mt} \omega^l ) -\sum_{i,j=1}^n \tilde{a}_{ij} \partial_{ij}(e^{-mt} \omega^l )  -\sum_{i=1}^n b_i \partial_i (e^{-mt} \omega^l ) \\
= & -m e^{-mt} \omega^l + l e^{-mt} \omega^{l-1} \omega_t \\
& -le^{-mt}\omega^{l-1}\sum_{i,j=1}^n \tilde{a}_{ij} \omega_{ij}   
 -l (l-1)e^{-mt} \omega^{l-2} \sum_{i,j=1}^n \tilde{a}_{ij} \omega_i \omega_j
-l e^{-mt}\omega^{l-1}\sum_{i=1}^n b_i \omega_i \\
= & e^{-mt} \omega^{l-2} \big( l\omega (\omega_t -\sum_{i,j=1}^n \tilde{a}_{ij} \omega_{ij} -\sum_{i=1}^n b_i \omega_i)  -l (l-1)\sum_{i,j=1}^n \tilde{a}_{ij} \omega_i \omega_j -m \omega^2 \big) 
\end{split}
\end{equation}

with
\begin{equation} \label{eq:omega:deriv}
\begin{split}
\omega_t (x,y,t+\tau_0) = & (18-\theta)\Lambda_t =(18-\theta)\Lambda \left( \frac{\theta}{(t+\tau_0)^2} -\frac{1}{t+\tau_0} \right), \\
\omega_i (x,y,t+\tau_0)= & -\theta_i \Lambda +(18-\theta) \frac{1}{4\pi (t+\tau_0)} \partial_i e^{-\frac{\theta}{t+\tau_0}}
= -\left( \frac{18-\theta}{t+\tau_0}+1 \right) \Lambda \theta_i ,\\
\omega_{ij} (x,y,t+\tau_0)= & -\left( \frac{18-\theta}{t+\tau_0}+1 \right) \Lambda \theta_{ij} 
+\frac{1}{t+\tau_0}\left( \frac{18-\theta}{t+\tau_0}+2 \right) \Lambda \theta_i \theta_j .
\end{split}
\end{equation}

Plugging \eqref{eq:omega:deriv} into the equation \eqref{eq:omega}, we have
\begin{equation}
\begin{split}
\mathcal{L} v = & l e^{-mt} \omega^{l-2} \Lambda^2 \Bigg[ (18-\theta)^2  \left( \frac{\theta}{(t+\tau_0)^2} -\frac{1}{t+\tau_0} +\frac{1}{t+\tau_0} \sum_{i,j=1}^n \tilde{a}_{ij}\theta_{ij} -\frac{1}{(t+\tau_0)^2} \sum_{i,j=1}^n \tilde{a}_{ij} \theta_i \theta_j +\frac{1}{t+\tau_0} \sum_{i=1}^n b_i \theta_i \right) \\
& -\frac{m}{l} (18-\theta)^2 +(18-\theta) \left( \sum_{i,j=1}^n \tilde{a}_{ij} \Bigg(  \theta_{ij} 
-\frac{2}{t+\tau_0} \theta_i \theta_j \Bigg) 
 +\sum_{i=1}^n b_i \theta_i \right) 
 - (l-1)\sum_{i,j=1}^n \tilde{a}_{ij} \left( \frac{18-\theta}{t+\tau_0}+1 \right)^2 \theta_i \theta_j   \Bigg]
\end{split}
\end{equation}

As \(l>1\), we see that
\begin{equation}
\mathcal{L} v \leq l e^{-mt} \omega^{l-2} \Lambda^2 [(18-\theta)^2  I +(18-\theta) II ]
\end{equation}
where
\begin{equation}
\begin{split}
I = & -\frac{\theta}{(t+\tau_0)^2} -\frac{1}{t+\tau_0} +\frac{1}{t+\tau_0} \sum_{i,j=1}^n \tilde{a}_{ij}\theta_{ij} -\frac{1}{(t+\tau_0)^2} \sum_{i,j=1}^n \tilde{a}_{ij} \theta_i \theta_j +\frac{1}{t+\tau_0} \sum_{i=1}^n b_i \theta_i -\frac{m}{l}, \\
II = & \sum_{i,j=1}^n \tilde{a}_{ij} \theta_{ij} 
-\frac{2l}{t+\tau_0} \sum_{i,j=1}^n \tilde{a}_{ij} \theta_i \theta_j 
 +\sum_{i=1}^n b_i \theta_i 
\end{split}
\end{equation}

For \( \theta=\overline{d}_\gamma^2=\frac{( x-x_0)^2 }{x+x_0} +\gamma^2\sum_{i=2}^n (y_i-(y_0)_i)^2 \), we have
\begin{equation} \label{eq:deriv:theta}
\begin{split}
\theta_1 = & \theta_x =\frac{2( x-x_0) }{x+x_0}-\frac{( x-x_0)^2 }{(x+x_0)^2}=\frac{( x+3x_0)( x-x_0) }{(x+x_0)^2}, \\
\theta_j = & \theta_{y_j} =2\gamma^2 (y_j-(y_0)_j), \ j=2,\ldots,n, \\
\theta_{11} = & \frac{( x-x_0) }{(x+x_0)^2} +\frac{( x+3x_0) }{(x+x_0)^2}-\frac{2( x+3x_0)( x-x_0) }{(x+x_0)^3} =\frac{8x_0^2 }{(x+x_0)^3} , \\
\theta_{ii} = & 2\gamma^2, \ i=2,\ldots,n, \\
\theta_{ij} = & 0, \ i\neq j \text{ and } \ i,j=1,\ldots,n.
\end{split}
\end{equation}

We observe that
\begin{equation}
|\theta|,|\theta_1|, |\theta_j|, |\theta_{11}|, |x\theta_{11}|, |\theta_{ii}| \leq C(\gamma) \text { on } K_{3\sqrt{2}}(P_0)
\end{equation}
when \(x_0, |y_0|\leq 1\). As \( \omega=(18-\theta) \frac{1}{4\pi t} e^{-\frac{\theta}{t}}\) and \(v(x,y,t)=e^{-mt} \omega^l(x,y,t+\tau_0)-M(\tau_0)\), we also have
\begin{equation}
\Vert v \Vert_{C^{1,1}} \leq C(\nu,\lambda).
\end{equation}

The conditions on the coefficients \eqref{assumption:coefficients} and the diagonality of \(D^2 \theta\) in \eqref{eq:deriv:theta} imply that
\begin{equation} \label{ineq:sum:deriv:theta}
\begin{split}
\sum_{i,j=1}^n \tilde{a}_{ij}\theta_{ij} =& \sum_{i=1}^n \tilde{a}_{ii}\theta_{ii} \leq \lambda^{-1} (x \theta_{11} +\sum_{i=2}^n \theta_{ii} ), \\
\sum_{i,j=1}^n \tilde{a}_{ij}\theta_i \theta_j =& \sum_{i=1}^n \tilde{a}_{ii}\theta_i^2 \geq \lambda ( x\theta_1^2 +\sum_{i=2} \theta_i^2), \\
|b_i|\leq & \lambda^{-1}, i=1,\ldots,n, \ \text{ and } |b_1| \geq \frac{\nu\lambda}{2}.
\end{split}
\end{equation}

Combining \eqref{eq:deriv:theta} and \eqref{ineq:sum:deriv:theta}, we can estimate the part \(I\) by
\begin{equation}
I\leq \frac{C(\gamma,\lambda,\nu)}{\tau_0^2 (\gamma)} -\frac{m}{l} \leq -\frac{m}{2l} 
\end{equation}
when we choose \(\frac{m}{l}\) sufficiently large, depending only on \( \gamma,\lambda\) and \(\nu\).

Also, the part \(II\) can be estimated by
\begin{equation}
II \leq \lambda^{-1} (x\theta_{11} + \sum_{i=2}^n \theta_{ii}+\sum_{i=2} |\theta_i| +\theta_1^+) - \frac{\nu\lambda}{2}\theta_1^- -c(\lambda, \gamma)l (x\theta_1^2+\sum_{i=2} |\theta_i|^2 )
\end{equation}

Like in the proof of Lemma 3.4 in \cite{daskalopoulos-lee03} in dimension 2, we can show that \(II\leq 0\) when \(d_\gamma \geq \frac{1}{4}\), and 
\begin{equation}
II \leq C(\nu, \lambda) \leq C(\nu, \lambda) (18-\theta)
\end{equation}
if  \(d_\gamma \leq \frac{1}{4}\), so we can see that
\begin{equation}
(18-\theta)^2  I +(18-\theta) II \leq  (18-\theta) \left(-\frac{m}{2l} +C(\nu, \lambda) \right) \leq 0
\end{equation}
with sufficiently large \(m\), concluding that \(\mathcal{L} v \leq 0\) in \(K_{3\sqrt{2}}\). We know that \(v\leq 0\) on \(\partial_p K_{3\sqrt{2}} \setminus \mathcal{Q}^1_{\frac{1}{2}} \) and \(v>c(\nu, \lambda)>0\) on \(\mathcal{Q}^2_{\frac{3}{2}}\). Setting \(\phi =\frac{v}{\inf_{ \mathcal{Q}^2_{\frac{3}{2}} }v}\), we see that \(\phi\) is the barrier function with \(\phi>1\) on \(\mathcal{Q}^2_{\frac{3}{2}}\) in the case \(\rho=1\).

Now we construct the barrier function \(\phi_\rho\) on \(\mathcal{Q}^2_{\frac{3\rho}{2}}\) for \(0<\rho<1\) by
\begin{equation}
\phi_\rho(\overline{d},t)=\phi \bigg( \frac{\overline{d}}{\rho},\frac{t}{\rho^2}\bigg).
\end{equation} 

Then \(\phi_\rho\) satisfies all the desired properties of the barrier function as 
\begin{equation}
\mathcal{L} \phi_\rho= (\phi_\rho)_t -L\phi_\rho=\frac{1}{\rho^2}( \phi_t -L \phi) \leq 0 \text { in } K_{3\sqrt{2}\rho} \setminus \mathcal{Q}^1_{\frac{\rho}{2}}
\end{equation}
and it satisfies the conditions \eqref{ineq1:phi-rho} as well as
\begin{equation} \label{ineq:C11:phi-rho}
\Vert \phi_\rho \Vert_{C^{1,1}(K_{3\sqrt{2}\rho} )} =\frac{1}{\rho^2} \Vert \phi\Vert_{C^{1,1}(K_{3\sqrt{2}} )} \leq \frac{1}{\rho^2} C(\nu,\lambda).
\end{equation}

\end{proof}

\subsection{Harnack inequality}

The following Harnack inequality is the \(n-\)dimensional version of Theorem 3.5 in \cite{daskalopoulos-lee03}. The quantity \(\rho_\nu \) is \eqref{nota:rho-nu}. We use the variable \(s=\sqrt{x}\) instead of \(x\) and \(s_0=\sqrt{x_0}\).

\begin{lemma} [Harnack inequality] \label{lemma:harnack}
Let \(u\geq 0\) be a classical solution of \(\mathcal{L}_s u =g\) in \(\mathcal{Q}_{\frac{\rho}{2}}(s_0, y_0, t_0 )\), where \(g\) is a bounded and continuous function in \(\mathcal{Q}_{\frac{\rho}{2}}(s_0, y_0, t_0 )\). Then it holds that
\begin{equation} \label{ineq:harnack}
\sup_{\mathcal{Q}_{\frac{\rho}{2}}(s_0, y_0, t_0 -\frac{3\rho^2}{4})} u \leq C\left( \inf_{\mathcal{Q}_{\frac{\rho}{2}}(s_0, y_0, t_0 )} u+\rho^{\frac{n}{n+1}}\rho_\nu(s_0)^{\frac{1}{n+1}} \Vert g \Vert_{L^{n+1} (\mathcal{Q}_{\frac{\rho}{2}}(s_0, y_0, t_0 ),d\mu)} \right)
\end{equation}
\end{lemma}

\begin{proof}
We follow along the lines of the proof of the elliptic Harnack inequality, Theorem 2.6 in \cite{daskalopoulos-lee03} with modifications about the dimension \(n\), and with the
A-B-P estimate, Lemma \ref{lem:ABP:lin}, and the barrier function from Lemma \ref{lem:barrier:ineq}. 
\end{proof}

In order to prove Lemma \ref{lemma:harnack}, we need the following lemma first. In the following lemma, we use the normalized measure \(|\mathcal{A}|_\mu \) of a set \(\mathcal{A}\) with respect to the measure \(d\mu =s^{\nu-1} ds dy dt\) 
\begin{equation}
|\mathcal{A}|_\mu=\frac{\gamma^{n}\nu}{2^{n}}\int_\mathcal{A} s^{\nu-1} ds dy dt.
\end{equation}

For example, the cube \(\mathcal{Q}_\rho(s_0,y_0,t_0)\) has the measure
\begin{equation}
\begin{split}
|\mathcal{Q}_\rho(s_0,y_0,t_0)|_\mu 
= & \frac{\gamma^{n}\nu}{2^{n}}\int_{t_0-\frac{\rho}{\gamma}}^{t_0+\frac{\rho}{\gamma}} \int_{(y_0)_n-\frac{\rho}{\gamma}}^{(y_0)_n+\frac{\rho}{\gamma}} \cdots \int_{(y_0)_2-\frac{\rho}{\gamma}}^{(y_0)_2+\frac{\rho}{\gamma}}  \int_{s_0-\rho_0}^{s_0+\rho} s^{\nu-1} ds dy_2 \
\ldots dy_n dt \\
= & [(s_0+\rho)^\nu -\overline{s_0}^\nu ]\rho^n
\end{split}
\end{equation}
where \(\overline{s_0}=\max(s_0-\rho,0)\).

\begin{lemma} \label{lemma:set:lowerbound:ratio}
Let \(u\) be a classical subsolution of the equation \(\mathcal{L}_s u \leq g\) in \(\mathcal{Q}_{3\sqrt{2} \rho}(s_0,y_0)\). Then there exist constants \(\varepsilon_0>0, 0<k<1\) and \(K>1\) such that if \(u\geq 0\) in \(\mathcal{Q}_{3\sqrt{2} \rho}(s_0,y_0)\), \(\inf_{\mathcal{Q}^2_{\frac{3\rho}{2} } (s_0,y_0) }u\leq 1\), and
\begin{equation} \label{assumption:g:small}
\rho^{\frac{n}{n+1}}\rho_\nu^{\frac{1}{n+1}}(s_0) \Vert g\Vert_{L^{n+1}(\mathcal{Q}_{3\sqrt{2} \rho}(s_0,y_0,t_0)d\mu)} \leq \varepsilon_0,
\end{equation}
then
\begin{equation} \label{ineq:less-than-set}
| \{u \leq K \} \cap \mathcal{Q}_\rho(s_0,y_0)|_\mu\geq k |\mathcal{Q}_\rho(s_0,y_0,t_0)|_\mu.
\end{equation}

\end{lemma}

\begin{proof}
As we know that the base point is \((s_0,y_0,t_0)\), we simply denote for \(r>0\), \(\mathcal{Q}_r=\mathcal{Q}_r (s_0,y_0,t_0) \) and \( \mathcal{B}_r=\mathcal{B}_r (s_0,y_0,t_0) \). For the barrier function \(\phi_\rho\) in Lemma \ref{lem:barrier:ineq}, we let \( w = u -2\phi_\rho\). Then
\begin{equation}
\mathcal{L}_s w \leq g -2\mathcal{L}_s \phi_\rho \leq g \ \text { in } K_{3\sqrt{2}\rho} \setminus \mathcal{Q}^1_{\frac{\rho}{2}}.
\end{equation}

Also, we have \( w \geq 0\) on \( \partial_p K_{3\sqrt{2}\rho}\) and \(\inf_{\mathcal{Q}^2_{\frac{3\rho}{2} } (x_0,y_0) } w \leq -1\) by the properties of \(\phi_\rho\) \eqref{ineq1:phi-rho}, thereby the A-B-P estimate \ref{lem:ABP:lin} implies that
\begin{equation}
1 \leq \inf_{ K_{3\sqrt{2}\rho} }  w^{-} \leq \inf_{ K_{3\sqrt{2}\rho} \setminus \mathcal{Q}^1_{\frac{\rho}{2}} }  w^{-} \leq C(n,\lambda,\nu) \rho^{n/(n+1)} \rho_\nu^{1/(n +1)} \left( \int_{\Gamma^-} ( g -2\mathcal{L}_s \phi_\rho  )^{n+1}(s,y,t) s^{\nu-1} ds dy dt \right)^{1/(n+1)}
\end{equation}
for \(\rho_\nu\) \eqref{nota:rho-nu}, the set \(\Gamma^- \)  \eqref{nota:gamma}, and the variable \( z =\frac{s^{2-\nu}}{2-\nu}\). Because \(\phi_\rho\) satisfies \eqref{ineq2:phi-rho} and \eqref{ineq:C11:phi-rho}, we have
\begin{equation}
1 \leq  C \rho^{n/(n+1)} \rho_\nu^{1/(n +1)} \Vert g \Vert_{L^{n+1}(\mathcal{Q}_{3\sqrt{2} \rho}(s_0,y_0, t_0),d\mu)} +\frac{C}{\rho^2} \rho^{n/(n+1)} \rho_\nu^{1/(n +1)} | \Gamma^- \cap \mathcal{Q}^1_{\frac{\rho}{2}}  |_\mu^{1/(n+1)}
\end{equation}

Choosing \(\varepsilon>0\) sufficiently small that \( C \rho^{n/(n+1)} \rho_\nu^{1/(n +1)} \Vert g \Vert_{L^{n+1}(\mathcal{Q}_{3\sqrt{2} \rho}(s_0,y_0, t_0),d\mu)} \leq \frac{1}{2}\), we have
\begin{equation}
\frac{1}{2} \leq C\rho^{\frac{-n-2}{n+1}} \rho_\nu^{\frac{1}{n +1}} | \Gamma^- \cap \mathcal{Q}^1_{\frac{\rho}{2}} |_\mu^{1/(n+1)}
\end{equation}

We have
\begin{equation}
\sup_{\mathcal{Q}_{3\sqrt{2} \rho}(s_0,y_0)} u \leq C \varepsilon_0
\end{equation}
because of Lemma \ref{lem:ABP:lin} and the assumption \eqref{assumption:g:small} on \(g\). So we have \(u \leq 1 \leq K\) on \(\mathcal{Q}_{3\sqrt{2} \rho}(s_0,y_0)\) if \(\varepsilon_0\) is sufficiently small, meaning that
\begin{equation}
c \rho^{n+2} \rho_\nu^{-1} \leq | \{ u \leq K\} \cap \mathcal{Q}^1_{\frac{\rho}{2}}  |_\mu.
\end{equation}

In order to prove \eqref{ineq:less-than-set}, it is sufficient, as \(\mathcal{Q}^1_{\frac{\rho}{2}}  \subset \mathcal{Q}_{\rho}  \), to show that
\begin{equation}
C \rho^{n+2} \rho_\nu^{-1} \geq | \{ u \leq K\} \cap \mathcal{Q}_{\rho}  |_\mu.
\end{equation}

We can see that
\begin{equation}
\delta(\rho) :=  \frac{\rho_\nu}{\rho^{n+2}} | \{ u \leq K\} \cap \mathcal{Q}_{\rho}  |_\mu =\frac{1}{\rho^2} [(s_0+\rho)^{2-\nu} -{s_0}^{2-\nu} ] \cdot[(s_0+\rho)^\nu -\overline{s_0}^\nu ]
\end{equation}
that satisfies 
\begin{equation}
\delta(\rho) \leq \frac{(3\rho)^\nu (3\rho)^{2-\nu}}{\rho^2} \leq C(\nu)
\end{equation}
if \(s_0\leq 2\rho\),
and
\begin{equation}
\delta(\rho) \leq \frac{1}{\rho^2} [(s_0+\rho)^{2-\nu} -{s_0}^{2-\nu} ] \cdot[(s_0+\rho)^\nu -(s_0-\rho)^\nu ] \leq C(\nu)s_0^{1-\nu} s_0^{\nu-1} \leq C(\nu)
\end{equation}
if \(s_0> 2\rho\), therefore finishing the proof of the lemma.

\end{proof}

Then we have the following lemma, from which lemma \ref{lemma:harnack} follows directly.
\begin{lemma} \label{lemma:harnack}
Let \(u \geq 0\) be a classical subsolution of the equation \(\mathcal{L}_s u \leq g\) in \(\mathcal{Q}_{3\sqrt{2} \rho}(s_0,y_0)\), where \(g\) is a bounded and continuous function in \(\mathcal{Q}_{3\sqrt{2} \rho}(s_0,y_0)\). Then there exist constants \(\varepsilon_0>0 \) and \(C>0\) depending only on \(\lambda\) and \(\nu\) such that if  \(\inf_{\mathcal{Q}^2_{\frac{3\rho}{2} } (s_0,y_0,t_0) }u\leq 1\) and
\begin{equation} \label{assumption:g:small:same}
\rho^{\frac{n}{n+1}}\rho_\nu^{\frac{1}{n+1}}(s_0) \Vert g\Vert_{L^{n+1}(\mathcal{Q}_{3\sqrt{2} \rho}(s_0,y_0,t_0), d\mu)} \leq \varepsilon_0,
\end{equation}
then
\begin{equation} \label{ineq:sup:upperbound}
\sup_{\mathcal{Q}^2_{\frac{3\rho}{2} } (s_0,y_0, t_0-\frac{3\rho^2}{4}) }u\leq C
\end{equation}
\end{lemma}

\begin{proof}
The function \(u\) satisfies the conditions of Lemma 2.11 in \cite{daskalopoulos-lee03}, in which \(l_i=\sigma K_0^{-\varepsilon/2} \theta^{-\varepsilon i/2} \) with \(\theta=\frac{K_0}{K_0-1}>\) for some \(K_0>1\), so we have an integer \(i_0\) depending only on universal constants such that
\begin{equation} \label{ineq:lemma:harnack}
\sum_{i\geq i_0}l_i \leq \frac{3}{2}.
\end{equation}

We want to show that 
\begin{equation} \label{claim:lemma:harnack}
\sup_{\mathcal{Q}^2_{\frac{3\rho}{2} } (s_0,y_0,t_0) }u\leq \theta^{i_0-1} K_0
\end{equation}
which proves the lemma \ref{lemma:harnack}. Let us prove the claim \eqref{claim:lemma:harnack} by contradiction. If the claim is not true, then there is a point \(P_{i_0}\) such that
\begin{equation}
P_{i_0} \in \mathcal{Q}^2_{\frac{3\rho}{2}} (s_0,y_0,t_0) \text{ and } u(P_{i_0}) \geq \theta^{i_0-1} K_0.
\end{equation}

Then, by the parabolic version of Lemma 2.11 in \cite{daskalopoulos-lee03}, there is a point \(P_{i_0+1}\) such that
\begin{equation}
P_{i_0 +1} \in \mathcal{Q}^2_{l_{i_0} \rho} (P_{i_0 }) \text{ and } u(P_{i_0 +1}) \geq \theta^{i_0} K_0.
\end{equation}

Repeating the process, we obtain a sequence of points \(\{ P_i, i\geq i_0\} \) satisfying
\begin{equation} \label{contradiction1:proof:lemma:proof}
P_{i +1} \in \mathcal{Q}^2_{l_{i} \rho} (P_{i }), \text{ and } u(P_{i +1}) \geq \theta^{i} K_0 \text{ for all } i\geq i_0.
\end{equation}
Then each \(P_i=(s_i, y_y,t_i) \) satisfies
\begin{equation} \label{contradiction2:proof:lemma:proof}
P_i \in \mathcal{Q}^2_{3\rho  } (s_0,y_0,t_0)
\end{equation} 
by the inequality \eqref{ineq:lemma:harnack} and
\begin{equation}
\begin{split}
|s_i-s_0| & \leq |s_{i_0}-s_0|+\sum_{k=i_0}^{i-1} |s_{k+1} -s_k| \leq \frac{3\rho}{2}+\sum_{k \leq i_0} l_k \rho \leq 3\rho \\
\gamma|y_i-s_0| & \leq \gamma|y_{i_0}-s_0| +\gamma\sum_{k=i_0}^{i-1} |y_{k+1} -y_k| \leq \frac{3\rho}{2}+\sum_{k \leq i_0} l_k \rho \leq 3\rho \\
|t_i-t_0| & \leq |t_{i_0}-t_0|+\sum_{k=i_0}^{i-1} |t_{k+1} -t_k| \leq \frac{3\rho}{2}+\sum_{k \leq i_0} l_k \rho \leq 3\rho .
\end{split}
\end{equation}
Because the sequence \(P_i=(s_i, y_y,t_i) \) satisfying both \eqref{contradiction1:proof:lemma:proof} and \eqref{contradiction2:proof:lemma:proof} contradicts the continuity of the function \(u\), we must have \eqref{claim:lemma:harnack}.

\end{proof}

\subsection{Hölder estimates}
We consider degenerate equations in the form
\begin{equation}
Lu-u_t =g
\end{equation}

What we need is the oscillation lemma, a parabolic version of Lemma 2.14 in \cite{daskalopoulos-lee03}.
\begin{lemma} \label{lem:oscillation}
Let \(u\) be a classical supersolution of equation \(L_s u =g\) in  \( \mathcal{Q}_{\rho }(s_0, y_0, t_0 )\) where \(g \) is a bounded continuous fuction. Then there exists universal constants \(0<\theta<1\) and \(C>0\) for which
\begin{equation} 
\mathop{osc}_{\mathcal{Q}_{\frac{\rho}{2}}(s_0, y_0, t_0 )} u \leq \theta \ \mathop{osc}_{\mathcal{Q}_{\rho}(s_0, y_0, t_0 )} u +C\rho^{\frac{n}{n+1}} \rho_{\nu} (s_0)^{\frac{1}{n+1}} \Vert g\Vert_{L^{n+1} ( \mathcal{Q}_{\rho}(s_0, y_0, t_0 ),d\mu)}.
\end{equation}
\end{lemma}

\begin{proof}
Following the line of the proof of Lemma 2.14 in \cite{daskalopoulos-lee03} with the parabolic Harnack inequality \eqref{ineq:harnack}, we get the desired result. 
\end{proof}

Using the oscillation lemma, Lemma \ref{lem:oscillation}, we can easily show the following parabolic version of Theorem 2.15 in \cite{daskalopoulos-lee03}.
\begin{lemma} \label{lem:hoelder:lin}
Let \(u\) be a classical supersolution of equation \(L_s u =g\) in  \( \mathcal{Q}_{\rho }(s_0, y_0, t_0 )\) where \(g \) is a bounded continuous fuction. Then there exist universal constants \(C>0\) and \(0<\alpha<\frac{n}{n+1}\) depending only on \(\lambda, \nu\) and \(n\), for which
\begin{equation} \label{eq:lem:oscillation}
\mathop{osc}_{\mathcal{Q}_{\rho}(s_0, y_0, t_0 )} u \leq 
C \rho^\alpha \Big( \rho_0^{-\alpha} \ \sup_{\mathcal{Q}_{\rho}(s_0, y_0, t_0 )} |u| 
+\rho_0^{\frac{n}{n+1} -\alpha} (s_0+\rho_0)^{\frac{1}{n+1}} \Vert g\Vert_{L^{n+1} ( \mathcal{Q}_{\rho}(s_0, y_0, t_0 ),d\mu)} \Big).
\end{equation}
\end{lemma}

\begin{proof}
Set \( \omega (\rho)=\mathop{osc}_{\mathcal{Q}_{\rho}(s_0, y_0, t_0 )} u \). Then Lemma \ref{lem:oscillation} shows that
\begin{equation}
\omega \Big(\frac{\rho}{2} \Big) \leq \theta \omega( \rho) +k(\rho)
\end{equation}
with a universal constant \(0<\theta<1\) and
\begin{equation}
k(\rho) =\rho^{\frac{n}{n+1} } (s_0+\rho_0)^{\frac{1}{n+1}} \Vert g\Vert_{L^{n+1} ( \mathcal{Q}_{\rho}(s_0, y_0, t_0 ),d\mu)}.
\end{equation}

Since both \(\omega\) and \(k\) are non-decreasing in \(\rho\), Lemma 8.23 in \cite{gilbarg-trudinger} implies \eqref{eq:lem:oscillation}.
\end{proof}

Now, we state the proof of the Hölder inequality, Theorem \ref{theorem:hoelder:lin}.
\begin{proof} (Proof of the Hölder inequality, Theorem \ref{theorem:hoelder:lin})
By controlling \(\sup_{\mathcal{Q}_{\rho}(s_0, y_0, t_0 )} |u| \) in \eqref{eq:lem:oscillation} with \( \Vert u \Vert_{C^0 (\mathcal{C}_1) }\), we will have Theorem \ref{theorem:hoelder:lin}.
\end{proof}


\section{Schauder estimates} \label{sec:shorttime}

We show Schauder estimates in the follwing two steps. First, we consider linear equations in general dimension. We modify the proof of the two-dimensional case, Theorem 4.1 in \cite{daskalopoulos-hamilton99}. 

Second, we study degenerate equations with variable coefficients, of which the two dimensional case is Theorem 7.1 in \cite{daskalopoulos-hamilton99}. The proof for degenerate equations combines the existence for linear equations and a standard perturbation argument, as done for Theorem II.1.1 in Daskalopoulos-Hamilton 98 \cite{daskalopoulos-hamilton98}. We will show the existence for fully nonlinear equations as well.

\subsection{Barriers and derivative estimates for the model degenerate equation}

Let us think about the following model degenerate equation. We consider the operator 
\begin{equation} \label{eq:model:operator}
L_0 f=f_t -(xf_{xx} +\sum_{i=2}^n f_{y_i y_i} +v f_x).
\end{equation}
\begin{lemma}
If \( f\) is smooth and satisfies the diffusion equation
\begin{equation} \label{eq:model}
f_t =x f_{xx} +\sum_{i=2}^n f_{y_i y_i} +v f_x +g.
\end{equation}
with transport velocity \(v>0\), on the box \(\mathcal{B}_r\)
\begin{equation}
\mathcal{B}_r= \{ 0 \leq x \leq r^2, -r\leq y_i \leq r, 1-r^2 \leq t \leq 1\},
\end{equation}
and if \( |f|\leq B\) on the box \( \mathcal{B}_r\), then for any \( \gamma <1\) we have
\begin{equation}
|f_x| \leq \frac{CB}{r^2} \ \text{ and } \ |f_{r_i}| \leq  \frac{CB}{r^2} \ \text{ for } \ i=2,\ldots,n
\end{equation}
on the box \( \mathcal{B}_{\gamma r}\).
\end{lemma}

We construct a barrier function \( \varphi\) as
\begin{equation}
\varphi =\frac{1}{t(x+at)} +\frac{1+t}{(1-x^2)^2}
 +\frac{1}{\big(x+b\sum_{i=2}^n (1-y_i)^2 \big) \sum_{i=2}^n (1- y_i)^2}
  +\frac{1}{\big(x+b\sum_{i=2}^n (1+y_i)^2 \big) \sum_{i=2}^n (1+ y_i)^2}
\end{equation}
which satisfies the barrier inequality
\begin{equation} \label{ineq:barrier}
\varphi_t > x \varphi_{xx} +\sum_{i=2}^n \varphi_{y_i y_i} +v \varphi_x -Cx\varphi^2 +c\varphi^{3/2}.
\end{equation}

\begin{proposition} \label{prop:barrier}
Given any \( v>0\), there is a \(b>0\) such that the function
\begin{equation}
\varphi=  \frac{1}{\big(x+b\sum_{i=2}^n y_i^2 \big) \sum_{i=2}^n y_i^2}
\end{equation}
satisfies the barrier inequality \eqref{ineq:barrier} on \( \{ 0< y_i <2, i=2,\ldots, n\} \).
\end{proposition}

\begin{proof}
Given \( v>0\), we want to find constants \( C<\infty, c>0, b>0\) satisfying
\begin{align*}
x\varphi_{xx} +\sum_{i=2}^n \varphi_{y_i y_i} +v \varphi_x +c\varphi^{3/2} < Cx\varphi^2.
\end{align*}

We have
\begin{align*}
\varphi_x =\frac{-1}{\big(x+b\sum_{i=2}^n y_i^2 \big)^2 \sum_{i=2}^n y_i^2} \ \text{ and } \ 
\varphi_{x x} =\frac{2}{\big(x+b\sum_{i=2}^n y_i^2 \big)^3 \sum_{i=2}^n y_i^2} \ \text{ . }
\end{align*}

In the tangential directions,
\begin{align*}
\varphi_{y_i} = & \frac{-2by_i}{\big(x+b\sum_{k=2}^n y_k^2 \big)^2 \sum_{k=2}^n y_k^2} 
+\frac{-2 y_i}{\big(x+b\sum_{k=2}^n y_k^2 \big) \big( \sum_{k=2}^n y_k^2 \big)^2 } \ \text{ , } \\
\varphi_{y_i y_i} = & \frac{8b^2 y_i^2 }{\big(x+b\sum_{k=2}^n y_k^2 \big)^3 \sum_{i=2}^n y_k^2} 
+\frac{8b y_i^2 }{\big(x+b\sum_{k=2}^n y_k^2 \big)^2 \big( \sum_{i=2}^n y_k^2\big)^2 } 
+\frac{-2b }{\big(x+b\sum_{k=2}^n y_k^2 \big)^2 \sum_{i=2}^n y_k^2}   \\
& +\frac{-2 }{\big(x+b\sum_{k=2}^n y_k^2 \big) \big( \sum_{i=2}^n y_k^2 \big)^2 } 
+\frac{8 y_i^2 }{\big(x+b\sum_{k=2}^n y_k^2 \big) \big( \sum_{i=2}^n y_k^2\big)^3 }   \ \text{ and } \\
\sum_{i=2}^n \varphi_{y_i y_i} = & \frac{8b^2 }{\big(x+b\sum_{k=2}^n y_k^2 \big)^3 } 
+\frac{(10-2n)b  }{\big(x+b\sum_{k=2}^n y_k^2 \big)^2 \sum_{i=2}^n y_k^2 } 
+\frac{10-2n }{\big(x+b\sum_{k=2}^n y_k^2 \big) \big( \sum_{i=2}^n y_k^2 \big)^2 }  \ \text{ . } \\
\end{align*}

Also
\begin{align*}
\varphi^{3/2} = \frac{1}{\big(x+b\sum_{i=2}^n y_i^2 \big)^{3/2} \big( \sum_{i=2}^n y_i^2 \big)^{3/2} }
\leq \frac{b^{-1/2}}{\big(x+b\sum_{i=2}^n y_i^2 \big) \big( \sum_{i=2}^n y_i^2 \big)^2 }
\end{align*}
and we need
\begin{align*}
& \frac{2x}{\big(x+b\sum_{i=2}^n y_i^2 \big)^3 \sum_{i=2}^n y_i^2}
+\frac{10-2n +c b^{-1/2}}{\big(x+b\sum_{k=2}^n y_k^2 \big) \big( \sum_{i=2}^n y_k^2 \big)^2 }
+\frac{(10-2n)b  }{\big(x+b\sum_{k=2}^n y_k^2 \big)^2 \sum_{i=2}^n y_k^2 } 
+\frac{8b^2 }{\big(x+b\sum_{k=2}^n y_k^2 \big)^3 } \\
\leq & \frac{v}{\big(x+b\sum_{k=2}^n y_k^2 \big)^2 \sum_{k=2}^n y_k^2}
+\frac{Cx }{\big(x+b\sum_{k=2}^n y_k^2 \big)^2 \big( \sum_{i=2}^n y_k^2 \big)^2 } \ \text{ . }
\end{align*}

Multiplying both sides by \( \big(x+b\sum_{k=2}^n y_k^2 \big)^3 \big( \sum_{i=2}^n y_k^2 \big)^2 \), we find that we need
\begin{equation} \label{ineq:barrier:cond}
\begin{split}
& 2x \sum_{i=2}^n y_k^2 
+(10-2n +c b^{-1/2}) \big(x+b\sum_{k=2}^n y_k^2 \big)^2
+(10-2n)b \big(x+b\sum_{k=2}^n y_k^2 \big) \sum_{i=2}^n y_k^2 
+8b^2 \big( \sum_{i=2}^n y_k^2 \big)^2 \\
< & v \big(x+b\sum_{k=2}^n y_k^2 \big) \sum_{i=2}^n y_k^2
+Cx\big(x+b\sum_{k=2}^n y_k^2 \big)
\end{split}
\end{equation}
and the inequality \eqref{ineq:barrier:cond} holds if we first choose the number \( b>0\) sufficiently small depending on \(v\) and then we pick the constant \( c>0\) sufficiently small depending on \(v \) and \(b\), and we  finally take the constant  \(C>0\) sufficiently large, depending on \( v\) and \( b\).
\end{proof}

\begin{corollary} \label{cor:barrier}
Given any \( v>0\), there is a \(b>0\) such that the function
\begin{equation} \label{eq:barrier}
\varphi=  \frac{1}{\big(x+b\sum_{i=2}^n (1-y_i)^2 \big) \sum_{i=2}^n (1-y_i)^2} 
+\frac{1}{\big(x+b\sum_{i=2}^n (1+y_i)^2 \big) \sum_{i=2}^n (1+y_i)^2} 
\end{equation}
satisfies the barrier inequality \eqref{ineq:barrier} on \( \{ 0< y_i <2, i=2,\ldots, n\} \).
\end{corollary}

\begin{proof}
The barrier inequality \eqref{ineq:barrier} is preserved under any translations and flips.
\end{proof}

The rescaled function
\begin{equation}
\widetilde{f}(x,y,t) =\frac{1}{B}f(r^2 x, r y,r^2 t)
\end{equation}
solves the same equation \eqref{eq:model} on the unit box \( \mathcal{B}_1\) and satisfies \( |\widetilde{f} | \leq 1\). So we may assume that \( f\) is defined on \( \mathcal{B}_1\) and satisfies \( |f| \leq 1\).

\begin{lemma}
If \( f\) is smooth and satisfies the diffusion equation
\[ f_t = x f_{xx} + \sum_{i=2} f_{y_i y_i} + v f_x\]
with transport velocity \( v>0\), and if \( |f|\leq 1\) on the box \( \mathcal{B}_1\), the for any \( \gamma<1\)
\begin{equation}
|f_x| \leq C \ \text{ and } \ |f_{y_i}| \leq C \ \text{ for } \ i = 2, \ldots,n
\end{equation}
on the box \( \mathcal{B}_\gamma\), for some constant \( C\) depending on \(v>0\) and \( \gamma<1\) but not depending on \(f\).
\end{lemma}

\begin{proof}
We use the quantity \( X = (A+f^2) f_x^2\) to estimate \( f_x\) where \(A\) is a constant to be set later.  \(f_x\) satisfies
\[ (f_x)_t = x f_{xxx} + \sum_{i=2} f_{ x y_i y_i} + (v +1) f_{xx} \] and \( X\) satisfies
\begin{align*}
X_x = & 2(A+f^2) f_x f_{xx} +2f f_x^3, \\
X_{y_i}= & 2(A+f^2) f_x f_{x y_i} +2f f_x^2 f_{y_i}, \\
X_{xx} = & 2(A+f^2) f_x f_{xxx} +2(A+f^2) f_{xx}^2 +10f f_x^2 f_{xx} +2f_x^4, \\
X_{y_i y_i}= & 2(A+f^2) f_x f_{x y_i y_i}+2(A+f^2) f_{x y_i}^2 +8f f_x f_{y_i} f_{x y_i} +2f f_x^2 f_{y_i y_i} +2 f_x^2 f_{y_i}^2 
\end{align*}
and hence
\begin{align*}
X_t = & 2(A+f^2) f_x (f_x)_t +2f f_x^2 f_t  \\
= & 2(A+f^2) f_x \big( x f_{xxx} + \sum_{i=2} f_{ x y_i y_i} 
+(v +1) f_{xx} \big) 
   +2f f_x^2 \big( x f_{xx} + \sum_{i=2} f_{y_i y_i} + v f_x \big) \\
= & x X_{xx} -x \big( 2(A+f^2) f_{xx}^2 +10f f_x^2 f_{xx} +2f_x^4 \big) \\
& +\sum_{i=2} X_{y_i y_i} -\sum_{i=2} \big( 2(A+f^2) f_{x y_i}^2 +8f f_x f_{y_i} f_{x y_i} +2f f_x^2 f_{y_i y_i} +2 f_x^2 f_{y_i}^2 \big) \\
& +(v +1) X_x - 2(v +1)f f_x^3  +2f f_x^2 x f_{xx} 
 +2f f_x^2 \sum_{i=2} f_{y_i y_i} 
 +2v f f_x^3 \\
= & x X_{xx} +\sum_{i=2} X_{y_i y_i} +(v +1) X_x  \\
& -2x (A+f^2) f_{xx}^2
 -8x f f_x^2 f_{xx}
 -2x f_x^4
 -2f f_x^3   \\ 
& -\sum_{i=2} 2(A+f^2) f_{x y_i}^2
 -\sum_{i=2} 8f f_x f_{y_i} f_{x y_i}
 -\sum_{i=2} 2 f_x^2 f_{y_i}^2  \\
\end{align*}

Assume that \( A \geq 8\). Then
\begin{align*}
& -2x (A+f^2) f_{xx}^2 \leq -(16+2f^2) x f_{xx}^2, \\
& -\sum_{i=2} 2(A+f^2) f_{x y_i}^2 \leq -\sum_{i=2} (16+2f^2) f_{x y_i}^2 , 
\end{align*} 
\begin{align*}
& -2x (A+f^2) f_{xx}^2
 -8x f f_x^2 f_{xx}
 -2x f_x^4
\leq -16 x ( f_{xx} +\frac{1}{4}f f_x^2)^2
 +x f^2 f_x^4
 -2x f_x^4 \leq  -x f_x^4, 
\end{align*} 
\begin{align*}
 -\sum_{i=2} 2(A+f^2) f_{x y_i}^2
 -\sum_{i=2} 8f f_x f_{y_i} f_{x y_i}
 -\sum_{i=2} 2 f_x^2 f_{y_i}^2  
\leq & -16\sum_{i=2} (f_{x y_i} +\frac{1}{4}f f_x f_{y_i})^2
 +\sum_{i=2} (f f_x f_{y_i})^2
 -\sum_{i=2} 2 f_x^2 f_{y_i}^2  \\
\leq &  -\sum_{i=2} f_x^2 f_{y_i}^2  
\end{align*} 
and hence
\begin{align*}
X_t \leq & x X_{xx} +\sum_{i=2} X_{y_i y_i} +(v +1) X_x  
 -x f_x^4
 +2f |f_x^3 |  
\end{align*}
and for \( \widetilde{X}=X/B\) with \( B=A^{5/2}\)
\begin{align*}
\widetilde{X}_t \leq & x \widetilde{X}_{xx} +\sum_{i=2} \widetilde{X}_{y_i y_i} +(v +1) \widetilde{X}_x  
 -C x \widetilde{X}^2
 +c \widetilde{X}^{3/2}
\end{align*}
for any \(C<\infty\) and \(C>0\), by making \(A\) sufficiently large. By the maximum principle, \(\widetilde{X} \leq \varphi\) for the barrier function \( \varphi\) \eqref{eq:barrier}. As \( v+1>0\), if the maximum occurs at the boundary, the term \((v +1) \widetilde{X}_x \) satisfies 
\((v +1) \widetilde{X}_x \leq (v +1) \varphi_x \) from \(  \widetilde{X}_x \leq \varphi_x\). The estimate of \( f_y\) can be obtained similarly by considering the quantity
\( Y=(A+f^2)f_y^2\) with sufficiently large constant \(A \geq 8\). Then \( Y\) satisfies
\begin{align*}
Y \leq xY_{xx} +\sum_{i=2} Y_{y_i y_i} +v Y_x -\frac{Y^2}{4A^2}
\end{align*}
and the estimate \(Y \leq B\) follows.
\end{proof}

\subsection{Schauder estimate of model degenerate equation}

\begin{theorem} [Schauder estimate]
For any \(v>0\), \( 0<\alpha<1\) and \(r<1\) there is a constant \(C>0\) so that
\begin{equation}
\Vert f \Vert_{C_s^{2+\alpha} (\mathcal{B}_r) } \leq C ( \Vert f \Vert_{C_s^{0} (\mathcal{B}_1) }
+\Vert L_0 f \Vert_{C_s^{\alpha} (\mathcal{B}_1) })
\end{equation}
for all \(C^\infty\) smooth functions f on \(\mathcal{B}_1\).

\end{theorem}

We begin with the following polynomial approximation, where the polynomial is given by the Taylor polynomial of degree 1 in \(x\) and \(t\) and of degree 2 in \(y=(y_2,\ldots,y_n)\). The maximum of a function \( f\) on the parabolic box \(\mathcal{B}_r \) is denoted by \(\Vert f \Vert_r\).

\begin{theorem} [Polynomial Approximation] \label{theorem:poly:approx1}
There exists a constant \(C\) such that for every smooth function \(f\) on the box \(\mathcal{B}_s \) we can choose a polynomial \(p\) of degree 1 in \(x\) and \(t\) and of degree 2 in \(y=(y_2,\ldots,y_n)\), which satisfies for every \(r\leq s\)
\begin{equation} \label{eq:error:Taylor}
\Vert f-p\Vert_r \leq C\left[ \left(\frac{r}{s}\right)^3 \Vert f\Vert_s +s^2 \Vert L_0 f\Vert_s \right].
\end{equation}
\end{theorem}

\begin{proof}
As in the proof of Theorem 5.2 in \cite{daskalopoulos-hamilton99}, we choose a bump function \(\psi\) on the set
\begin{align*}
\mathcal{S}=\{ 0\leq x <\infty, -\infty<y_i<\infty, i=2,\ldots,n, -\infty<t\leq 1\}
\end{align*}
such that \( \psi=0\) outside \( \mathcal{B}_1\) and \( \psi=1\) on \( \mathcal{B}_{1/2}\). Let \(h\) be the unique smooth bounded solution on \(\mathcal{S}\) of the equation
\begin{align*}
L_0 h=\psi L_o f.
\end{align*}
For the error \(k=f-h\), we let \(p\) be the Taylor polynomial of \(k\) at the point \((x,y,t)=(0,\vec{0},1)\) given by
\begin{align*}
p(x,y,t)=k(0,\vec{0},1)+k_x(0,\vec{0},1)\cdot x+\sum_{i=2}^n k_{y_i} (0,\vec{0},1)\cdot y_i +k_t(0,\vec{0},1)\cdot(t-1)+\frac{1}{2}\sum_{i=2}^n k_{y_i y_i }(0,\vec{0},1) \cdot y_i^2
\end{align*}

Then the polynomial \(p\) satisfies the estimate \eqref{eq:error:Taylor}, if we follow the proof of Theorem I.6.1 in \cite{daskalopoulos-hamilton98} with a few obvious changes: (1) the operator \(L_0\) is replaced by \eqref{eq:model:operator}; (2) instead of a single variable \(y=y_2\) we use \((n-1)-\)dimentional vector variable \(y=(y_2,\ldots,y_n)\); (3) the polynomial \(p\) is not of degree one in \(y=y_2\) as in in \cite{daskalopoulos-hamilton98}, but of degree two in \(y=(y_2,\ldots,y_n)\) as in \cite{daskalopoulos-hamilton99}. \end{proof}

Then we have the follwing result in \cite{daskalopoulos-hamilton99} consequently.

\begin{theorem} \label{theorem:poly:approx2}
For each \(v>0\) and \(0<\alpha<1\), there exists a constant \(C\) with the following property. If \( f\) is a smooth function on the box \(\mathcal{B}_1\) whose Taylor polynomial at \((0,\vec{0},1)\) of degree 1 in \(x,t\) and 2 in \(y\) is zero, then
\begin{equation}
\sup_{0<r\leq 1} \frac{\Vert f\Vert_r}{r^{2+\alpha} } \leq C \left( \Vert f\Vert_1+ \sup_{0<r\leq 1} \frac{\Vert L_0 f\Vert_r}{r^{\alpha} }\right).
\end{equation}
\end{theorem}

\begin{proof}
Exactly as in the proof of Theorem 5.3 in \cite{daskalopoulos-hamilton99} and Theorem I.7.1 in \cite{daskalopoulos-hamilton98} with the obvious changes.
\end{proof}

For a smooth function \(f\) on the box \( \mathcal{B}_1\), let \(T_{12,1} f\) be the Taylor polynomial at \((0,\vec{0},1)\) of degree 1 in \(x,t\) and 2 in \(y\)
\begin{align*}
T_{12,1} f =f(0,\vec{0},1)+f_x(0,\vec{0},1)\cdot x+\sum_{i=2}^n f_{y_i} (0,\vec{0},1)\cdot f_i +f_t(0,\vec{0},1)\cdot(t-1)+\frac{1}{2}\sum_{i=2}^n f_{y_i y_i }(0,\vec{0},1) \cdot y_i^2
\end{align*}
and \(R_{12,1} f\) be the remainder
\begin{align*}
R_{12,1} f = f -T_{12,1} f .
\end{align*}
Also, for a smooth function \(g\) on the box \( \mathcal{B}_1\), let \( T_0 g\) denote the Taylor polynomial of degree 0 in both space and time at the point \( (0,\vec{0},1)\)
\begin{align*}
T_0 g =g (0,\vec{0},1),
\end{align*}
i.e. the constant given by its evaluation at the point \( (0,\vec{0},1)\), and let \(R_0 g\) be its remainder
\begin{align*}
R_0 g= g-T_0 g.
\end{align*}

\begin{corollary}
There exists a constant \(C\) such that for any smooth function \(f\) on the box \( \mathcal{B}_1\), 
\begin{equation}
\sup_{0<r\leq 1} \frac{\Vert R_{12,1} f\Vert_r}{r^{2+\alpha} } \leq C \left( \Vert R_{12,1} f\Vert_1+ \sup_{0<r\leq 1} \frac{\Vert R_0 L_0 f\Vert_r}{r^{\alpha} }\right).
\end{equation}
\end{corollary}

\begin{proof}
We apply Theorem \ref{theorem:poly:approx2} to the remainder \(R_{12,1} f\) where \(T_{12,1} R_{12,1} f=0\).
\end{proof}

We define a metric \( ds^2 :=\frac{d x^2 }{x} +\sum_{i=2}^n d y_i^2 \) with distance 
\[ s[(x_1, \vec{y_1}, t_1),(x_2, \vec{y_2}, t_2)] = s[(x_1, \vec{y_1}),(x_2, \vec{y_2})] +\sqrt{|t_1-t_2|}\] 
where
\begin{align*}
c(|\sqrt{x_1}-\sqrt{x_2}|+|\vec{y_1}-\vec{y_2}|)
\leq s[(x_1, \vec{y_1}),(x_2, \vec{y_2})]
\leq C(|\sqrt{x_1}-\sqrt{x_2}|+|\vec{y_1}-\vec{y_2}|).
\end{align*}

Using similar arguments as in \cite{daskalopoulos-hamilton99}, we have the following interior Schauder estimate in higher dimensions, which corresponds to Theorem 5.11 in \cite{daskalopoulos-hamilton99} in dimension two.

\begin{theorem} [Interior Schauder Estimate]
There exists \(\lambda >0 \) such that for any \(\mu <\lambda\) and any smooth function \(f\) on the parabolic cylinder \( \mathcal{C}_\lambda (Q) \) 
\begin{align*}
\Vert u\Vert_{\mathcal{C}_s^{2+\alpha} (\mathcal{C}_\mu(Q))} \leq C( \ \Vert u\Vert_{\mathcal{C}^{0} (\mathcal{C}_\lambda(Q))} +\Vert L_0 u\Vert_{H_s^{\alpha} (\mathcal{C}_\lambda (Q))} \ )
\end{align*}
\end{theorem}

By applying dilations and standard rescalings, and by differentiating the model equation multiple times, we obtain the Schauder estimate Theorem \ref{theorem:schauder}, which is Theorem 5.25 in \cite{daskalopoulos-hamilton99} in dimension two. In addition, we get the short-time existence as a corollary at the end.

\subsection{Short time existence of a solution to the model degenerate equation}

We recall the smoothing operator in \cite{daskalopoulos-hamilton99} first. Using the smoothing operator, we show the short time existence of smooth solution of the model equation.

Let \(\mathcal{S}_0\) be the half space of \( \mathbb{R}^n\) where \(x_1 \geq 0\) and \(\mathcal{S} =\mathcal{S}_0 \times [0,\infty)\). For \(T>0\), we let \(\mathcal{S}_T =\mathcal{S}_0 \times [0,T]\). For a point \(P=(x, \vec{y} ) \) on the half space \(S_0\), and \(Q=(u, \vec{v})\) any point in the unit box \(\mathcal{B}_1 =\{ |u|<1, \Vert \vec{v} \Vert  \leq 1 \} \), we define the point 
\begin{align*}
M_\varepsilon (P;Q)=(\xi, \vec{\zeta})=\big( ( \sqrt{x+2\varepsilon} +\sqrt{\varepsilon}u)^2, \vec{y} +\sqrt{\varepsilon} \vec{v} \big).
\end{align*}
Then we can see that \( s[ (x+2\varepsilon,\vec{y}), (\xi, \vec{y}) ] =\sqrt{\varepsilon} |u|, \ s[ (\xi,\vec{y}), (\xi, \vec{\zeta}) ] =\sqrt{\varepsilon} |v|\). 

Next, we let \(\varphi\) be a standard nonnegative bump function with support in \(\mathcal{B}_1\) with \(\int \varphi =1\). For a function \( h\) defined on the half space \(S_0\) the space regularization \(h_\varepsilon\) of \(h\) is defined as
\begin{align*}
h_\varepsilon(P) =\int_{Q =(u,\vec{v})\in \mathcal{B}_1} \varphi (u, \vec{v}) h(M_\varepsilon(P,((u, \vec{v})) ) du dv_2 \ldots dv_n
\end{align*}
for \(P=(x,\vec{y}) \in S_0\).

Following the regularization and extension lemmas as in \cite{daskalopoulos-hamilton99} with obvious changes with respect to the dimension, we get the following existence and uniqueness result, Theorem 6.4 in \cite{daskalopoulos-hamilton99}.

\begin{theorem} [Existence and uniqueness] \label{theorem:short:model}
Let \(k\) be a nonnegative integer and \(0<\alpha<1\). If \(g \in \mathcal{C}_s^{k,\alpha}(\mathcal{S})\) and \(f^0 \in \mathcal{C}_s^{k,2+\alpha}(\mathcal{S})\) with compact support in \(\mathcal{S}\) and \(\mathcal{S}_0\) respectively, then for any constant \(c\) and any \(v>0, T>0\), the initial value problem
\begin{align*}
\begin{cases}
L_0 f -cf =g  \text{ in } \mathcal{S}_T \\
f(\cdot,0)=f^0 \text{ on } S_0
\end{cases}
\end{align*}
admits a unique solution \(f \in \mathcal{C}_s^{k,2+\alpha}(\mathcal{S}_T) \) satisfying the estimate
\begin{align*}
\Vert f \Vert_{ \mathcal{C}_s^{k,2+\alpha} (\mathcal{S}_T)}
\leq C(T) (\Vert f^0 \Vert_{\mathcal{C}_s^{k,2+\alpha}(\mathcal{S})} +\Vert g \Vert_{\mathcal{C}_s^{k,\alpha}(\mathcal{S})} )
\end{align*}
for some constant \(C(T)\) which depends only on \(k,\alpha,v,c\) and \(T\).

\end{theorem}

\subsection{Variable coefficient degenerate equations}
Let \(\mathcal{D}=\mathcal{D}_1\) be the unit disk in \( \mathbb{R}^n\) and \(Q=\mathcal{D} \times [0,T]\) for some \(T>0\). We denote \(\mathcal{D}_+ =\mathcal{D} \cap \{x_1\geq 0\}\). We study the short time existence of a smooth solution to the following degenerate equation. 
\begin{equation} \label{eq1:model:variable}
w_t=\sum_{i,j=1}^n a^{ij} w_{ij} +\sum_{i=1}^n b^i w_i +cw
\end{equation}
where the coefficients \(a^{ij}, b^i, c\) belong to the Hölder class \(\mathcal{C}^{k,\alpha}(\mathcal{D}_{1-\delta/2} \times [0,T])\) for some constant \(0<\alpha<1\) and a nonnegative integer \(k\), for any \(0<\delta<1\). For such a \(\delta\), we have a collection of charts \(\gamma_l:\mathcal{D}_+ \rightarrow \mathcal{D}_\delta (P_l) \cap \mathcal{D} \) which flatten the boundary of \( \mathcal{D} \), where \(P_l\) are finite numbers of points on the the boundary of \( \mathcal{D} \) such that \(P_l=\gamma_l(0)\) and \( \mathcal{D}_\delta (P_l) \) is the disk of radius \(\mathcal{D}_\delta (P_l)\) with center \(P_l\), such that the collection of points \( \mathcal{D}_\delta (P_l) \) covers the boundary of \( \mathcal{D} \).

We assume that there is a number \(\delta\) so that for every \(l \in I\), each coordinate change \(\gamma_l\) transforms 
the operator \( L\) on \(  \mathcal{D}_\delta (P_l) \cap \mathcal{D} \)
\begin{equation} \label{operator1:model:variable}
L[w]:=w_t- (x_1 a^{11} w_{11} +2\sqrt{x_1} \sum_{i=2}^{n} a^{1i} w_{1i} +\sum_{i,j=2}^n a^{ij} w_{ij} +\sum_{i=1}^n b^i w_i )
\end{equation}
into an operator \( \widetilde{L}_l \) on \( \mathcal{D}_+\)
\begin{equation} \label{operator2:model:variable}
\widetilde{L}_l[\widetilde{w}]:=\widetilde{w}_t- (x_1 \widetilde{a}^{11} \widetilde{w}_{11} +2\sqrt{x_1} \sum_{i=2}^{n} \widetilde{a}^{1i} \widetilde{w}_{1i} +\sum_{i,j=2}^n \widetilde{a}^{ij} \widetilde{w}_{ij} +\sum_{i=1}^n \widetilde{b}^i \widetilde{w}_i )
\end{equation}
where the coefficients \( \widetilde{a}^{ij}, \widetilde{b}^i\) belong to the Hölder class \(\mathcal{C}_s^{k,\alpha}\) and they satisfy the conditions
\( \sum_{i,j=1}^n \widetilde{a}^{ij} \xi_i \xi_j \geq \lambda |\xi|^2 \) for any \(\xi \in \mathbb{R}^n\), \(|\widetilde{b}_{i} |\leq \lambda^{-1}\) for \(i=1,\ldots,n\), and \(\widetilde{b}_{1} \geq \lambda\)  for some number \( \lambda >0\). 

By applying a perturbation argument from Theorem II.1.1 in \cite{daskalopoulos-hamilton98} to Theorem  \ref{theorem:short:model} with obvious changes in dimension \(n\geq 3\), we get the following existence and uniqueness result which corresponds to the two-dimension version Theorem 7.1 in \cite{daskalopoulos-hamilton99}.

\begin{theorem} [Existence and uniqueness] \label{theorem:short:variable}
Let \(k\) be a nonnegative integer, \(0<\alpha<1\) a number, and \( L\) the operator \eqref{operator1:model:variable}  satisfying the properties aove on \(\mathcal{Q} \). If \(g \in \mathcal{C}_s^{k,\alpha}(\mathcal{Q})\) and \(f^0 \in \mathcal{C}_s^{k,2+\alpha}(\mathcal{D})\), then for any constant \(c\) and any \(v>0, T>0\), the initial value problem for the operator \(L\) 
\begin{align*}
\begin{cases}
L f =g  \text{ in } Q \\
f(\cdot,0)=f^0 \text{ on } \Omega
\end{cases}
\end{align*}
admits a unique solution \(f \in \mathcal{C}_s^{k,2+\alpha}(\mathcal{Q}_T) \) satisfying the estimate
\begin{align*}
\Vert f \Vert_{ \mathcal{C}_s^{k,2+\alpha} (\mathcal{Q})}
\leq C(T) (\Vert f^0 \Vert_{\mathcal{C}_s^{k,2+\alpha}(\mathcal{D})} +\Vert g \Vert_{\mathcal{C}_s^{k,\alpha}(\mathcal{Q})} )
\end{align*}
for some constant \(C(T)\) which depends only on \(k,\alpha, \lambda\) and \(T\).
\end{theorem}

\section*{\bf{Acknowledgements}}
Ki-Ahm Lee  was supported by the National Research Foundation of Korea (NRF) under the grant NRF-
2021R1A4A1027378. Ki-Ahm Lee also holds a joint appointment with Research Institute of Mathematics of Seoul National University.

\end{document}